\documentstyle[amssymb,12pt]{article}
\def\rom#1{{\normalshape #1}}
\newtheorem{proposition}{Proposition}
\newtheorem{conjecture}{Conjecture}
\newenvironment{proof}{\begin{trivlist}\item[]
{\bf Proof.} }{\hfill $\Box$ \end{trivlist}}
\begin{document}
\begin{center}{\LARGE\bf Cayley Hypersurfaces}\\[20pt]
{\large\bf Michael Eastwood and Vladimir Ezhov}\end{center}
{\small {\bf Abstract:} We exhibit a family of homogeneous hypersurfaces in
affine space, one in each dimension, generalising the Cayley surface.}

The Cayley surface in affine three space is given by
$$\textstyle x_3=x_1x_2-\frac{1}{3}x_1{}^3.$$
See, for example,~\cite[Chapter~III \S6]{ns} for a discussion of its
properties. In $N$-dimensions, we may consider the hypersurface given by
\begin{equation}\label{cayley}
\Phi_N(x_1,x_2,\ldots,x_N)\equiv\sum_{d=1}^N(-1)^d{\textstyle\frac{1}{d}}
   \sum_{\makebox[0pt]{\makebox[20pt][l]{$\scriptstyle i+j+\cdots+m=N$}}}
     \overbrace{x_ix_j\cdots x_m}^d=0.\end{equation}
This is the Cayley surface when~$N=3$. The next few are as follows.
$$\begin{array}l
x_4=x_1x_3+\frac{1}{2}x_2{}^2-x_1{}^2x_2+\frac{1}{4}x_1{}^4\\[5pt]
x_5=x_1x_4+x_2x_3-x_1{}^2x_3-x_1x_2{}^2+x_1{}^3x_2-\frac{1}{5}x_1{}^5\\[5pt]
x_6=x_1x_5+x_2x_4+\frac{1}{2}x_3{}^2-x_1{}^2x_4-2x_1x_2x_3-\frac{1}{3}x_2{}^3
    +x_1{}^3x_3+\frac{3}{2}x_1{}^2x_2{}^2-x_1{}^4x_2+\frac{1}{6}x_1{}^6.
\end{array}$$
Since the first term in (\ref{cayley}) is $-x_N$ and this is the only
occurrence of this variable, these hypersurfaces are polynomial graphs over the
remaining variables.

The Cayley surface is affine homogeneous. This follows immediately from
$\Phi_3$ being annihilated by the following two linearly independent affine
vector fields:
$$\frac{\partial}{\partial x_1}+x_1\frac{\partial}{\partial x_2}+
x_2\frac{\partial}{\partial x_3}\quad\mbox{and}\quad
\frac{\partial}{\partial x_2}+x_1\frac{\partial}{\partial x_3}.$$
The hypersurface defined by (\ref{cayley}) generalises sufficiently many
properties of the Cayley surface that we call it the Cayley hypersurface.
That it is affine homogeneous is an immediate consequence of the following:
\begin{proposition} The polynomial $\Phi_N(x_1,x_2,\ldots,x_N)$ is annihilated
by the vector fields
\begin{equation}\label{fields}X_p\equiv\frac{\partial}{\partial x_p}+
           \sum_{h=p+1}^Nx_{h-p}\frac{\partial}{\partial x_h}\end{equation}
for $p=1,2,\ldots,N-1$.\end{proposition}
\begin{proof} We compute
$$\frac{\partial}{\partial x_h}\Phi_N=\sum_{d=1}^N(-1)^d
\sum_{\makebox[0pt]{\makebox[20pt][l]{$\scriptstyle j+\cdots+m=N-h$}}}
     \overbrace{x_j\cdots x_m}^{d-1}\quad\mbox{for }h=1,2,\ldots N$$
with the convention that a product with no terms is $1$ (so that when $h=N$,
this formula gives~$-1$ and, otherwise, we can start the sum at~$d=2$).
Therefore,
$$\begin{array}{rcl}X_p\Phi_N&=&\displaystyle\sum_{d=2}^N(-1)^d
\sum_{\makebox[0pt]{\makebox[20pt][l]{$\scriptstyle j+\cdots+m=N-p$}}}
     \overbrace{x_j\cdots x_m}^{d-1}
 \;+\sum_{h=p+1}^Nx_{h-p}\sum_{d=1}^N(-1)^d
\sum_{\makebox[0pt]{\makebox[20pt][l]{$\scriptstyle j+\cdots+m=N-h$}}}
     \overbrace{x_j\cdots x_m}^{d-1}\\[22pt]
&=&\displaystyle\sum_{d=2}^N(-1)^d
\sum_{\makebox[0pt]{\makebox[20pt][l]{$\scriptstyle j+\cdots+m=N-p$}}}
     \overbrace{x_j\cdots x_m}^{d-1}
 \;+\;\sum_{d=1}^N(-1)^d
\sum_{\makebox[0pt]{\makebox[20pt][l]{$\scriptstyle j+\cdots+m=N-p$}}}
     \overbrace{x_ix_j\cdots x_m}^{d}.
\end{array}$$
These expressions evidently cancel. \end{proof}
\begin{proposition} The Cayley hypersurface admits a transitive Abelian group
of affine motions.\end{proposition}
\begin{proof} The vector fields (\ref{fields}) commute and so may be
exponentiated to the required Abelian group. (In fact, this is how the Cayley
surface is defined in~\cite[p.~93]{ns}.) \end{proof}
\begin{proposition} The isotropy algebra of the Cayley hypersurface is
generated by
\begin{equation}\label{H}
H\equiv\sum_{h=1}^Nhx_h\frac{\partial}{\partial x_h}.\end{equation}
\end{proposition}
\begin{proof} Each term in $\Phi_N$ is weighted homogeneous of weight $N$ if
$x_h$ has weight~$h$. It follows that $H\Phi_N=N\Phi_N$ and, in particular,
$H\Phi_N|_{\{\Phi_N=0\}}=0$. We are required to prove that, up to scale, $H$ is
the only vector field of the form
$$X=\sum_{i=1}^N\sum_{j=1}^Na_{ij}x_i\frac{\partial}{\partial x_j}$$
with this property. Since $\Phi_N$ is an irreducible polynomial of degree~$N$
and $X\Phi_N$ is a polynomial of degree at most~$N$, if $X\Phi_N$
vanishes along $\{\Phi_N=0\}$, then $X\Phi=c\Phi_N$ for some constant~$c$.
Therefore, after adding a suitable multiple of $H$, it suffices to show that
if $X\Phi_N=0$, then~$X=0$. Notice that $\Phi_N$ has the following form
$$\pm\underbrace{x_1{}^N/N}_{\deg=N}
\mp\underbrace{\vphantom{/}x_1{}^{N-2}x_2}_{\deg=N-1}
\pm\underbrace{\vphantom{/}x_1{}^{N-3}x_3+p(x_1,x_2)}_{\deg=N-2}
\mp\underbrace{\vphantom{/}x_1{}^{N-4}x_4+q(x_1,x_2,x_3)}_{\deg=N-3}
\pm\cdots$$
for suitable polynomials $p(x_1,x_2)$, $q(x_1,x_2,x_3)$, \ldots. By
considering firstly the leading term $x_N{}^N/N$, it follows that $X$ cannot
have any terms in $\partial/\partial x_1$. Alternatively, it is this
consideration which determines which multiple of $H$ to use in the initial
modification. Then, by looking at the term of degree $N-1$, we see that $X$
cannot involve $\partial/\partial x_2$. Then the terms of degree $N-2$ dispense
with $\partial/\partial x_3$ and so on. Working through $\Phi_N$ in this way,
it follows that~$X=0$. \end{proof}
\begin{proposition} The affine normals of the Cayley hypersurface are
everywhere parallel to the $x_N$-axis.
\end{proposition}
\begin{proof} At the origin we have
$$x_N=\sum_{i,j}g^{ij}x_ix_j+\sum_{i,j,k}a^{ijk}x_ix_jx_k+\cdots$$
where
\begin{equation}\label{ga}
g^{ij}=\cases{1&if $i+j=N$\cr 0&else}\quad\mbox{ and }\quad
a^{ijk}=\cases{1&if $i+j+k=N$\cr 0&else.}\end{equation}
Thus $g_{ij}$, the inverse of $g^{ij}$, is given by the same formula and
$\sum_{i,j}g_{ij}a^{ijk}=0$. In fact, all the higher order tensors are
trace-free too. This property of the cubic terms characterises the affine
normal~\cite{l}. Note from (\ref{fields}) that the symmetries $X_p$
all commute with $\partial/\partial x_N$. It follows that they preserve the
$x_N$-direction but, on the other hand, since the affine normal is affinely
invariant, these symmetries take one affine normal to another. For the isotropy
symmetry (\ref{H}) we have $[\partial/\partial x_N,H]=N(\partial/\partial x_N)$
which says that the corresponding $1$-parameter subgroup simply rescales $x_N$
along its axis. \end{proof}

We conjecture these various properties are enough to characterise these
hypersurfaces:
\begin{conjecture}
Suppose $\Sigma$ is a non-degenerate hypersurface in affine $N$-space
such\linebreak that:
\vspace*{-7pt}
\begin{itemize}\addtolength{\itemsep}{-7pt}
\item $\Sigma$ admits a transitive Abelian group of affine motions
\item The full symmetry group of $\Sigma$ has one-dimensional isotropy
\item The affine normals to $\Sigma$ are everywhere parallel\newline
      \rom{(}i.e.~$\Sigma$ is an `improper affine hypersphere'\/\rom{)}.
\end{itemize}
\vspace*{-7pt}
Then $\Sigma$ is given by \rom{(\ref{cayley})} in a suitable affine
co\"ordinate system.
\end{conjecture}

We have verified this for hypersurfaces in four dimensions (as a special case
of classifying the homogeneous non-degenerate hypersurfaces with isotropy or
classifying those with a transitive Abelian group of affine motions). Details
will appear elsewhere. If two-dimensional isotropy is allowed, then another
variation on the Cayley surface arises, namely
$$\textstyle x_4=x_1x_3+\frac{1}{2}x_2{}^2-\frac{1}{3}x_1{}^3.$$
This sort of variation is discussed in~\cite{dv} and~\cite[pp.~121--122]{ns}.

In~\cite{np}, Nomizu and Pinkall give a differential geometric characterisation
of the Cayley surface: it is not assumed {\em a priori\/} that the surface is
homogeneous. It is not clear how to extend this characterisation to Cayley
hypersurfaces.

Yet another generalisation of Cayley surface, is suggested by Dillen
and Vrancken~\cite{DV} as a hypersurface with parallel difference tensor
together with some genericity condition. The explicit defining equation (6.3)
of \cite{DV} very much resembles~(\ref{cayley}). Curiously, there is a whole
family of homogeneous hypersurfaces
$$\sum_{d=1}^N(-1)^d{\textstyle\frac{1}{d!}}
  \prod_{n=0}^{d-3}\left[(1-b)n+2\right]
   \sum_{\makebox[0pt]{\makebox[20pt][l]{$\scriptstyle i+j+\cdots+m=N$}}}
     \overbrace{x_ix_j\cdots x_m}^d=0$$
interpolating between them. If $b=0$ this is (\ref{cayley}) and if $b=1$ it is
(6.3) of~\cite{DV}.

The Cayley surfaces are ruled. This property also extends to hypersurfaces:
\begin{proposition} If $N$ is odd, then the Cayley hypersurface is uniquely
ruled by $(N-1)/2$-planes. If $N$ is even, then it is uniquely ruled by
$(N-2)/2$-planes.\end{proposition}
\begin{proof} If $N$ is odd and we fix $x_1,x_2,\ldots,x_{(N-1)/2}$, then
(\ref{cayley}) is linear in the remaining variables. If $N$ is even and we fix
$x_1,x_2,\ldots,x_{N/2}$, then (\ref{cayley}) is linear in the remaining
variables. In both cases, uniqueness follows by examining the quadratic terms
which are non-degenerate of split signature. These determine the possible
directions in which a maximal embedded plane may point, only one of which is
consistent with the higher order terms.\end{proof}

For surfaces, the Pick invariant is precisely the third order obstruction to
its being ruled. Though there is no such obstruction in higher dimensions, we
have:
\begin{proposition} The Cayley hypersurfaces have vanishing Pick invariant.
\end{proposition}
\begin{proof} With reference to~(\ref{ga}), the Pick invariant is
$$\sum_{i,j,k,l,m,n}g_{il}g_{jm}g_{kn}a^{ijk}a^{lmn}$$
which evidently vanishes.\end{proof}

{\bf Acknowledgement.} We would like to thank Franki Dillen, Udo Simon, and
Takeshi Sasaki for interesting comments and encouragement.

\small\renewcommand{\section}{\subsection}

\sc
\begin{flushleft}
Department of Pure Mathematics\\
University of Adelaide\\
South AUSTRALIA 5005\\[5pt]
\normalshape
E-mail: meastwoo@spam.maths.adelaide.edu.au\\
\mbox{\phantom{E-mail: }vezhov@spam.maths.adelaide.edu.au}
\end{flushleft}
\end{document}